\documentclass[a4paper,leqno]{article}

\usepackage{amsmath}
\usepackage{amsfonts}
\usepackage{enumerate}
\usepackage{theorem}

\theoremstyle{plain}
\newtheorem{teo}{Theorem}
\newtheorem{lem}[teo]{Lemma}

\newtheorem{prop32}{Proposition}[section]

{\theorembodyfont{\rmfamily}

}

\renewcommand{\eqref}[1]{\textnormal{(\ref{#1})}}

\numberwithin{equation}{section}

\newcommand{\cvd}{\hfill$\square$}
\newcommand{\proof}[1]{\noindent\textsc{Proof#1}}

\title{Corrigendum to ``Determining a sound-soft polyhedral scatterer
by a single far-field measurement''}

\author{\begin{tabular}{lcl}
Giovanni Alessandrini & and & Luca Rondi\\
\normalsize{\texttt{alessang@units.it}} &
& \normalsize{\texttt{rondi@units.it}}
\end{tabular}\\
\normalsize{Dipartimento di Matematica e Informatica}\\
\normalsize{Universit\`a degli Studi di Trieste, Italy}}

\date{}

\begin{document}

\setcounter{section}{3}
\setcounter{secnumdepth}{1}
\setcounter{prop32}{1}

\maketitle

In the paper, \cite{Ale e Ron}, on the determination of a sound-soft polyhedral scatterer by a single far-field
measurement, the proof of Proposition 3.2 is incomplete. In this corrigendum we provide a new proof of the same
proposition which fills the previous gap. In order to introduce it, we recall some definitions from \cite{Ale e Ron}.

Let $v$ be a nontrivial real
valued solution to the Helmholtz equation
\begin{equation}\label{eq}
\Delta v + k^2v=0 \ \text{in }G ,
\end{equation}
in a connected open set $G\subset \mathbb{R}^N$, $N\geq 2$.
We denote the \emph{nodal set} of $v$
as
$$\mathcal{N}_v=\{x\in G:\ v(x)=0\}$$
and we let $\mathcal{C}_v$ be the set of \emph{nodal critical
points}, that is
$$\mathcal{C}_v=\{x\in G:\ v(x)=0\text{ and }\nabla v(x)=0\}.$$
We say that $\Sigma\subset \mathcal{N}_v$ is a \emph{regular
portion} of $\mathcal{N}_v$ if it is an analytic open and
connected hypersurface contained in
$\mathcal{N}_v\backslash\mathcal{C}_v$. Let us denote by
$A_1,A_2,\ldots,A_n,\ldots$ the \emph{nodal domains} of $v$ in $G$,
that is the connected components of $\{x\in G:\ v(x)\neq
0\}=G\backslash \mathcal{N}_v$.
Let us recall the statement of Proposition 3.2 in \cite{Ale e Ron}.

\begin{prop32}[\cite{Ale e Ron}]\label{orderprop}
We can order the nodal domains $A_1,A_2,\ldots,A_n,\ldots$ in such
a way that for any $j\geq 2$ there exist $i$, $1\leq i<j$, and a
regular portion $\Sigma_j$ of $\mathcal{N}_v$ such that
\begin{equation}\label{order}
\Sigma_j\subset \partial A_i\cap\partial A_j.
\end{equation}
\end{prop32}

The gap in the proof given in \cite{Ale e Ron} stands in the fact that the ordering $A_1,A_2,\ldots,A_n,\ldots$
obtained with that method might not ensure that all the nodal domains are contained in the sequence.
We base the new proof on the following theorem.

\begin{teo}
The set $\mathcal{C}_v$ has Hausdorff dimension not exceeding $N-2$.
\end{teo}

A proof can be found in \cite[Theorem~2.1]{Lin}. Further developments of the theory on the structure of zero sets of solutions
to elliptic equations can be found, for instance, in \cite{Han e Har e Lin,Har e Hof e Nad} and in their references.

Let $G'=G\backslash\mathcal{C}_v$. By the property of $\mathcal{C}_v$ described in the previous theorem,
and by using \cite[Chapter~VII, Section~4]{Hur e Wal} and \cite[Theorem~IV 4, Corollary~2]{Hur e Wal},
we can conclude that $G'$ is an open and connected set.
We also remark that, for every $x\in\mathcal{N}_v\backslash\mathcal{C}_v$,
there are exactly two nodal domains, $A$ and $B$, of $v$ such that $x\in\partial A\cap\partial B$.
Finally, let us note that the nodal domains of $v$ in $G$ coincide with the nodal domains of $v$ in $G'$.

We shall also make use of the following elementary lemma.

\begin{lem}\label{lemma1}
For any connected open set $G\subset\mathbb{R}^N$, there exists an increasing sequence
$\{G_m\}_{m=1}^{\infty}$ of bounded, connected open sets such that $G=\bigcup_{m=1}^{\infty}G_m$ and
$G_m\subset\subset G$ for every $m$.
\end{lem}

\proof{.} For every $k=1,2,\ldots$, we denote
$$D_k=\{x\in G:\ \mathrm{dist}(x,\partial G)>1/k,\ |x|<k\}.$$
Let us assume, without loss of generality, that $D_1\neq \emptyset$ and let us fix $y\in D_1$.
For every $x\in \overline{D_k}$, let $\gamma_x$ be a path in $G$ joining $y$ to $x$.
For every $h>0$, let $\mathcal{U}^h_x=\{z\in\mathbb{R}^N:\ \mathrm{dist}(z,\gamma_x)<h\}$.
We obviously have that $\mathcal{U}^h_x$ is a connected open set. Let $h(x)>0$ be such that
$\mathcal{U}^{h(x)}_x\subset\subset G$. We have that $\{\mathcal{U}^{h(x)}_x\}_{x\in \overline{D_k}}$
is an open covering of the compact set $\overline{D_k}$.
Therefore, we can find $x_1,\ldots,x_l\in \overline{D_k}$ such that
$\overline{D_k}\subset\bigcup_{j=1}^l\mathcal{U}^{h(x_j)}_{x_j}$.
We observe that $E_k=\bigcup_{j=1}^l\mathcal{U}^{h(x_j)}_{x_j}$ is an open connected set such that
$\overline{D_k}\subset E_k\subset\subset G$. Therefore the lemma follows choosing
$G_m=\bigcup_{k=1}^m E_k$.\cvd

\bigskip

\proof{ of Proposition~\ref{orderprop}.}
We apply Lemma~\ref{lemma1} to the connected set $G'=G\backslash\mathcal{C}_v$.
We choose $A_1$ such that $A_1\cap G_1\neq \emptyset$ and we proceed by induction.

Let us assume that we have ordered $A_1,\ldots,A_n$ in such a way
that there exist $\Sigma_2,\ldots,\Sigma_n$ regular portions of
$\mathcal{N}_v$ such that \eqref{order} holds for any
$j=2,\ldots,n$ and for some $i<j$.

Let $\hat{A}_n=\stackrel{\circ}{\overline{A_1\cup\ldots\cup A_n}}$. If
$G'\backslash \hat{A}_n=\emptyset$, then we are done. Otherwise, let $m\geq 1$ be the smallest number such that
$G_m\backslash \hat{A}_n\neq \emptyset$. Since $G_m$ is connected, 
we can find $y\in\partial \hat{A}_n\cap G_m$ and $r>0$ such that $B_r(y)\cap\partial \hat{A}_n$ is
a regular portion of $\mathcal{N}_v$ and there exist exactly two nodal domains, $\tilde{A}_1\subset \hat{A}_n$
and $\tilde{A}_2$ with $\tilde{A}_2\cap \hat{A}_n=\emptyset$, whose intersections with $B_r(y)$ are not empty.
Clearly, $\tilde{A}_1$ coincides with $A_i$, for some $i=1,\ldots,n$, and if we pick $A_{n+1}=\tilde{A}_2$
and $\Sigma_{n+1}=B_r(y)\cap \mathcal{N}_v$, then \eqref{order} holds for $j=n+1$, too.

If $G$ contains only finitely many nodal domains, then we can iterate this construction and after a finite number of steps
we recover all the nodal domains, that is for some $l\in\mathbb{N}$ we have $G'\backslash \hat{A}_l=\emptyset$
and we are done.
Otherwise, we argue in the following way.
Since $\overline{G_m}$ is contained in $G'$, for every $x\in \overline{G_m}$ there is a neighbourhood of
$x$ intersecting at most two different nodal domains. By compactness, we obtain that
$\overline{G_m}$ intersects at most finitely many different nodal domains. Hence, if we iterate the previous construction,
after a finite number of steps we find $l\in\mathbb{N}$ such that
$G_m\backslash \hat{A}_l=\emptyset$. By repeating the argument for the smallest $m'>m$ such that
$G_{m'}\backslash \hat{A}_l\neq\emptyset$, we conclude that for any $m\in\mathbb{N}$
there exists $l\in\mathbb{N}$ such that
$G_m\backslash \hat{A}_l=\emptyset$. Therefore the infinite sequence $\{A_i\}$ comprises all the nodal domains
of $v$ in $G$.\cvd

\subsection{Acknowledgements}
The authors wish to express their gratitude to Hongyu Liu and Jun Zou for pointing out to them the gap in the proof
of Proposition~3.2 in \cite{Ale e Ron} and for kindly sending them their preprint \cite{Liu e Zou}.

\end{document}